\documentclass[draft]{tbimc}
\newtheorem{thm}{Theorem}[section]

\newtheorem{nas}{Corollary}[section]
\theoremstyle{remark}
\newtheorem{zau}{Remark}[section]
\newtheorem{example}{Example}[section] 
\theoremstyle{definition}
\newtheorem{ozn}{Definition}[section]

\DeclareMathOperator{\pr}{\mathsf P} % !!!!!=======
\DeclareMathOperator{\M}{\mathsf E}
\begin{document}

\selectlanguage{english}  %%%%% ===========  !!!!!!! ПРОСТАВЬ В КОНЦЕ !!!!!

\title[WAVELET-BASED SIMULATION OF RANDOM PROCESSES]{WAVELET-BASED SIMULATION\\
 OF RANDOM PROCESSES FROM CERTAIN CLASSES\\ 
       WITH GIVEN ACCURACY AND RELIABILITY}

\author{Ievgen Turchyn}
\address{Department of Mechanics and Mathematics, 
Oles Honchar Dnipro National University, 
 Gagarin av., 72, Dnipro,  49010, Ukraine}
\curraddr{}
\email{evgturchyn@gmail.com}

%    General info
\UDC{519.21}
\subjclass[2010]{Primary 60G10; Secondary 42C40}
\date{??.??.2019}
\dedicatory{}
\keywords{Wavelets, Sub-Gaussian random processes, Simulation}
%% \thanks{}

\begin{abstract}
   We consider stochastic processes $Y(t)$ which can be represented 
as $Y(t)=(X(t))^s, s \in \mathbb{N},$ 
where $X(t)$ is a stationary  strictly  sub-Gaussian process
and  build a wavelet-based model that simulates $Y(t)$
 with given accuracy and reliability  in $L_p([0,T])$. 
 A model for simulation 
 with given accuracy and reliability  in $L_p([0,T])$
is also built for processes $Z(t)$ which can be represented 
as $Z(t)=X_1(t) X_2(t)$, where $X_1(t)$ and 
$X_2(t)$ are  independent stationary strictly  sub-Gaussian processes. 
\end{abstract}

\maketitle

%%%%%%%%%%%%%%%%%%%%%%%%%%%%%%%%%%%%%%%%%%%%%%%%

%\selectlanguage{english}
%% \setcounter{page}{74} 

%\title[Рівномірна збіжність розкладів для $\varphi$-субгауссових процесів]
    %{ Умови рівномірної збіжності розкладів $\varphi$-субгауссових випадкових процесів по системам функцій, породженим вейвлетами }

%%    Information for author
%\author{Ю. В. Козаченко}
%\address{ Київський національний університет ім. Тараса Шевченко,
%вул. Володимирська,  64, 
%Київ, 01033, Україна}
%\email{yvk@univ.kiev.ua}
%
%%    Information for second author
%\author{Є. В. Турчин}
%\address{Кафедра вищої математики,
    %факультет механізації сільського господарства,
    %OLES Honchar DNU,    
    %Dnipro,
    %Україна}  
%%\curraddr{(????) %Київський державний університет технологій і дизайну,
%%вул.~Немировича-Дан\-ченко, 2,02011, Київ
%%    }
%\email{evgturchyn@gmail.com}
%
%
%%    General info
%\subjclass[2010]{Primary 60G10; Secondary 42C40} %%% ВСЁ ОК 2010 
%\date{??/??/20??}
%%\dedicatory{Цю статтю присвячено нашим вчителям.}
%\keywords{Wavelets, Sub-Gaussian random processes, Simulation}
%
%\begin{abstract}
   %В роботі побудовано розклади з некоррельованими    
%коефіціентами випадкових процесів другого порядку 
%по системам функцій, що породжуються вейвлетами. Для розкладів, 
%коефіціенти яких незалежні, знайдено умови їх рівномірної
%збіжності з імовірністю 1 на скінченому інтервалі. Для розкладів строго 
%$\varphi$-субгауссових процесів знайдено умови рівномірної
%збіжності за ймовірністю на скінченому інтервалі.    
%\end{abstract}
%
%\maketitle
%
\section{Introduction}

  Wavelet expansions and  wavelet-based expansions form an interesting 
class of representations of random processes. At present there exist 
many articles devoted to such expansions and their properties, some of them are 
\cite{Taqqu-rateopt}, \cite{Camb-Masry-94}, \cite{Did-Pipiras},
\cite{Koz-Olenko-Polosm-SAA-2011}, \cite{Koz-Olenko-Polosm-CommSt-2014} and 
\cite{koz-turchyn-IJSMS}--\cite{Meyer-FBM}.    	
Wavelet-based expansions with uncorrelated terms (see, for instance, 
articles \cite{Did-Pipiras}, \cite{Meyer-FBM}) are especially important 
since they are very convenient for approximation and 
simulation of random processes.  

   We will consider simulation of stochastic processes with given 
accuracy and reliability. This is simulation of a random process by 
a model which has guaranteed rate of convergence in a certain sense --- 
i.e., a model $\hat{X}(t)$ approximates a process $X(t)$   	 
with given accuracy $\varepsilon$ and reliability	$1-\delta$
in a functional space $U(\mathbf{T})$
if 
$$
P\{\|X - \hat{X}\|_{U(\mathbf{T})} > \varepsilon\}  \leq \delta.  
$$

  Simulation of stochastic processes with given 
accuracy and reliability has been studied,
 in particular, in \cite{mono-Elsev}
(see also e.g. \cite{Koz-Dovgay-Roz}, \cite{Koz-Pogor}, \cite{Koz-Sottin} 
and \cite{koz-turchyn-IJSMS}).  
It is necessary to mention that results on simulation
with given accuracy and reliability are available mostly for 
light-tailed processes  --- Gaussian and sub-Gaussian processes
(although there are some exceptions, see, for instance, \cite{Turchyn-ModSt}).  

	The article is devoted to simulation with given accuracy
and reliability in $L_p([0,T])$ of random processes which 
can be represented as 
$$
 Y(t) = (X(t))^s, \ s \in \mathbb{N}, 
$$ 
where $X(t)$ is a stationary strictly sub-Gaussian  process, and 
$$
   Z(t) = X_1(t) X_2(t), 
$$
where $X_1(t)$ and $X_2(t)$ are stationary strictly sub-Gaussian  processes.  
So our approach allows simulation of processes with one-dimensional 
distributions which have relatively heavy tails
(i.e. heavier than the Gaussian ones). 

  Our models are derived from a model which was studied in \cite{koz-turchyn-IJSMS}, 
where a wavelet-based expansion was considered and used for 
construction of a model of a process 
(this or similar 
wavelet-based expansions and their rate of convergence 
were also studied in \cite{Did-Pipiras}, 
\cite{Walt-Zh} and \cite{Zhao}).    			

%%   MENTION [Turchyn-ModStoch] !!!

   %Питання рівномірної збіжності випадкових рядів з імовірністю 1 та 
%за ймовірністю розглядалось в  роботах \cite{B-K}, \cite{Kahane}, 
%\cite{Koz-Perest-Vas}, \cite{Walter-Shen} та інших. 
%
   %В цій роботі отримані розклади $\varphi$-субгауссових випадкових процесів 
%з некоррельованими коефіціентами по системам функцій, що породжуються вейвлетами. Знайдено умови рівномірної збіжності з імовірністю 1 на скінченому інтервалі для розкладів з незалежними коефіціентами і рівномірної збіжності
%за ймовірністю на скінченому інтервалі для розкладів строго 
%$\varphi$-субгауссових процесів. 
%
   %Робота складається з 6 розділів. В другому розділі наведено основні відомості
%з теорії строго $\varphi$-субгауссових випадкових величин та процесів. Третій
%розділ містить теорему про розклад центрованого випадкового процесу другого
%порядку в ряди з некоррельованими коефіціентами по системам функцій, 
%побудованих за допомогою систем вейвлетів. В четвертому розділі наведені умови рівномірної збіжності з імовірністю 1 отриманих розкладів в ряди з 
%незалежними членами на скінченому інтервалі. В п'ятому розділі знайдено 
%умови рівномірної збіжності за ймовірністю на скінченому інтервалі розкладів
%строго $\varphi$-субгауссових випадкових процесів.
%
    %Всі одержані в роботі результати є новими.
%

\section{Sub-Gaussian random variables and processes} 
% РОЗДІЛ 2 
% Add \boldmath ?!

%\begin{ozn} \cite{Kr-Rut}
  %Неперервна парна опукла функція $\varphi=\{ \varphi(x), x \in R \}$  
  %називається $N$-функцією Орліча, якщо $\varphi(0)=0$, 
%$\varphi(x)>0$ коли $x \neq 0$
%і виконуються умови:
%$$
 %\lim\limits_{x \rightarrow 0} \frac{\varphi(x)}{x}=0,\,
 %\lim\limits_{x \rightarrow \infty} \frac{\varphi(x)}{x}=\infty.
%$$
%\end{ozn}
%%
%\begin{ozn} \cite{Gi-K-N} 
   %Будемо говорити, що для $N$-функції  $\varphi$ виконується умова Q, якщо 
%$$
%\lim \limits_{x \rightarrow 0}\inf 
%\frac{\textstyle {\varphi(x)}}{\textstyle {x^2}}=c>0
%$$
%(зокрема, $c$ може дорівнювати $+\infty$).
%\end{ozn}

\begin{ozn} \cite{B-K}
   Let $\{ \Omega,\mathfrak{F} , \pr \}$ be a standard probability space. 
% стандартний імовірнісний простір.
 A random variable $\xi$ is called sub-Gaussian if 

1) $\M\xi=0$; 

2) for all $\lambda \in \mathbb{R}$ there  exists  $\M \exp\{\lambda \xi \}$; 

3) there exists such constant $a>0$ that inequality
$$
\M \exp\{\lambda \xi \} \leq \exp\{ a^2 \lambda^2/2  \} 
$$
holds for all  $\lambda \in \mathbb{R}$.
\end{ozn}

We will denote the set of all sub-Gaussian random variables
by ${\rm Sub}(\Omega)$.
   ${\rm Sub}(\Omega)$ is a Banach space with respect to the norm 
$$
   \tau(\xi)=\sup\limits_{\lambda \neq 0}
\left(\frac{ 2\ln \M \exp\{\lambda \xi \} }{\lambda^2}\right)^{1/2}
$$
(see \cite{B-K}).
   Examples of sub-Gaussian random variables can be found in 
 \cite{B-K}. Let us note that centered normal random variables belong 
to ${\rm Sub}(\Omega)$.

\begin{ozn} \cite{B-K}
A sub-Gaussian random variable $\xi$ is called strictly sub-Gaussian
if $\tau (\xi)=(\mathsf{E} \xi^2)^{1/2}$.
\end{ozn}

\begin{ozn} \cite{B-K}
A family of sub-Gaussian random variables $\Delta$ is called 
strictly sub-Gaussian if for any finite or countable set 
$\{\xi_i, i \in I \}$ from $\Delta$ and all $\lambda_{i} \in 
\mathbb{R}$ holds relation 
\[
\tau^2 \left(\sum\limits_{i \in I} \lambda_{i} \xi_i
 \right)=\mathsf{E} \left(\sum\limits_{i \in I} \lambda_{i} \xi_i
 \right)^2 .
\]
\end{ozn}

\begin{ozn} \cite{B-K}
A stochastic process $X=\{X(t), t \in \mathbf{T} \}$ 
is called strictly sub-Gaussian if the family of random variables
$\{ X(t), t \in \mathbf{T} \}$ is strictly sub-Gaussian.
\end{ozn}

\begin{example} 
\cite{B-K}
 Let $X=\{X(t), t \in \mathbf{T} \}$ be a centered Gaussian
process. Then $X(t)$ is a strictly sub-Gaussian stochastic process.
\end{example}

\begin{example} %% \cite{Koz-Pashko}
 Let $X=\{X(t), t \in \mathbf{T} \}$ be a random process 
such that 
$$
X(t)=\sum\limits_{k=1}^{\infty} \xi_k f_{k}(t), 
$$
 where
$\xi=\{\xi_k, k=1,2,\ldots \}$ is a family of independent 
strictly sub-Gaussian random variables and for all 
$t \in \mathbf{T}$
\[
\sum\limits_{k=1}^{\infty} \mathsf{E} \xi_k^2 f_{k}^2 (t) < \infty.
\] 
Then $X(t)$ is a strictly sub-Gaussian stochastic process.
\end{example}

%
%
%\begin{ozn} \cite{Koz-Kov} Сім'я $\Delta$ випадкових величин
%$\xi \in Sub_{\varphi}(\Omega)$ називається строго $\varphi$-субгауссовою,
%якщо існує константа $C_{\Delta}>0$ така, що для будь-якої скінченої множини
%$I$ випадкових величин $\xi_i \in \Delta$ та довільних $\lambda_i \in R$ має
%місце нерівність:
%$$
%\tau^2_{\varphi} \left( \sum_{i \in I} \lambda_i \xi_i  \right) \leq
%C_{\Delta}^2 \M \left( \sum_{i \in I} \lambda_i \xi_i \right)^2. 
%$$
   %Константа $C_{\Delta}$ називається визначальною константою сім'ї $\Delta$.
%\end{ozn}
%
%\begin{lem} \cite{Koz-Kov} Лінійне замикання строго $\varphi$-субгауссової
%сім'ї $\Delta$ в просторі $L_2(\Omega)$ є строго $\varphi$-субгауссовою 
%сім'єю з тією самою визначальною константою.
%\end{lem}
%
%\begin{ozn} \cite{B-K} Випадковий процес $X=\{X(t), t \in T \}$ називається 
%$\varphi$-субгаус- \linebreak совим, якщо всі випадкові величини $X(t), t \in T,$ є
%$\varphi$-субгауссовими.
%\end{ozn}
%
%\begin{ozn} \cite{Koz-Kov} Випадковий процес $X=\{X(t),\,\, t \in T \}$ 
%називається строго \linebreak 
%$\varphi$-субгаус\-совим, якщо сім'я випадкових величин 
%$\{X(t), t \in T \}$ є строго $\varphi$-суб\-гаус\-со\-вою з визначальною константою
%$C_T$.  Константа $C_T$ називається визначальною константою процесу $X$.
%\end{ozn}
%
  
% Examples of sub-Gaussian and strictly sub-Gaussian random variables 
% and random processes can be found in \cite{B-K}. 

\section{Expansion of a random process into a wavelet-based series}
% РОЗДІЛ 3  

\begin{ozn} 
\cite{hardle1998}
Let $\phi \in L_2(\mathbb{R})$ be such a function that the
following assumptions hold:

i)
\[
\sum\limits_{k \in \mathbb{Z}} |\hat \phi(y+2 \pi k)|^2=1
\]
almost everywhere, where $\hat{\phi}(y)$ is the Fourier transform of
$\phi$;
%\[
%\hat \phi(y)=\int\limits_{\mathbb{R}} \exp\{-iyx\} \phi(x) dx ;
%\]

ii) There exists a function $m_0 \in L_2([0, 2\pi])$ such that $m_0(x)$
has period $2\pi$ and almost everywhere
\[
\hat \phi(y)=m_0\left(\frac{y}{2}\right)
 \hat\phi \left(\frac{y}{2}\right) ;
\]

iii)  $\hat \phi(0) \not= 0$ and the function $\hat \phi(y)$
is continuous at 0.

  The function $\phi(x)$ is called $f$-wavelet. Let $\psi(x)$ be the inverse
Fourier transform of the function
\[
\hat \psi(y)=\overline{m_0\left(\frac{y}{2}+\pi \right)}
\exp\Bigl \{\Bigr. -i\frac{y}{2} \Bigl\}\Bigr.
\hat \phi \left(\frac{y}{2} \right) .
\]
The function $\psi(x)$ is called $m$-wavelet.
\end{ozn}

 Let
\[
\phi_{jk}(x)=2^{j/2}\phi(2^{j}x-k), \,\,\,
\psi_{jk}(x)=2^{j/2}\psi(2^{j}x-k), \,\,\, k\in \mathbb{Z},
\,\,\, j=0,1,2,\ldots
\]
The family of functions $\{ \phi_{0k}, \psi_{jk},
j=0,1,\ldots, k \in \mathbb{Z} \}$ is an orthonormal basis in
 $L_2(\mathbb{R})$
 (see,
for example,  \cite{hardle1998}).

 \begin{zau}
    We will consider only real-valued wavelets below. 		
 \end{zau} 

%Any function $f \in L_2(\mathbb{R})$ can be represented in the form
%\begin{equation}
%\label{1.1}
%f(x)=\sum\limits_{k \in \mathbb{Z}} \alpha_{0k}\phi_{0k}(x) +
%\sum\limits_{j=0}^{\infty}\sum\limits_{k \in \mathbb{Z}} \beta_{jk}\psi_{jk}(x),
%\end{equation}
%where $\alpha_{0k}=\int\limits_{\mathbb{R}}f(x)\overline{\phi_{0k}(x)}dx$,
%$\beta_{jk}=\int\limits_{\mathbb{R}}f(x)\overline{\psi_{jk}(x)}dx$,
%\[
%\sum\limits_{k \in \mathbb{Z}}|\alpha_{0k}|^2+
%\sum\limits_{j=0}^{\infty}\sum\limits_{k \in \mathbb{Z}}|\beta_{jk}|^2 < \infty.
%\]
%That is, series (\ref{1.1}) converges in the norm of the space
 %$L_2(\mathbb{R})$.
%Representation (\ref{1.1}) is called wavelet representation.

   Let us now formulate a result which is 
very important for us. 		

\begin{thm}  \cite{Koz-Turch-Roz}
%%%%% ================== W-exp-n, general process ========= 
\label{p3-cor-rozkl1-and-3} % НАСЛІДОК ПРО РОЗКЛАДИ № 1 І 3
    Suppose that $X(t),\, t \in \mathbb{R}$, is a  centered second-order
random process such that its correlation function
$R(t,s)=\mathsf{E} X(t) \overline{X(s)}$ can be re\-pre\-sen\-ted as
\begin{equation}  
\label{COR-F-INT-REPR}
R(t,s) =
\int_{\mathbb{R}}u(t, y)\overline{u(s,y)}dy,
\end{equation}
where $u(t,y)$  is a Borel function which belongs to $L_{2}(\mathbb{R})$ for all   
$t\in \mathbb{R}, \{ \phi_{0k}(x), \linebreak \psi_{jk}(x),\,\,$
  $k \in \mathbb{Z},\, j=0,1, \ldots \}$  is an arbitrary wavelet basis,
\begin{equation} % 
\label{p3eq-rozkl1-a0k-gen}
a_{0k}(t)=
\frac{1}{\sqrt{2\pi}} \int_{\mathbb{R}} u(t,y) \overline{\hat \phi_{0k}(y)}dy, 
%=
%\frac{1}{\sqrt{2\pi}} \int_{\mathbb{R}}
% u(t,y) \overline{\hat \phi(y)}\exp\{iyk\}dy,
\end{equation}
\begin{equation}
\label{p3eq-rozkl1-bjk-gen}
b_{jk}(t)=\frac{1}{\sqrt{2\pi}}
\int_{\mathbb{R}} u(t,y) \overline{\hat \psi_{jk}(y)}dy,
%=
%\frac{1}{\sqrt{2\pi}2^{j/2}} \int_{\mathbb{R}} u(t,y)\exp\left\{i\frac{y}{2^j}k \right\}
%\overline{\hat \psi(y/2^j)}dy,
\end{equation}
$\hat \phi_{0k}(y)$ and $\hat \psi_{jk}(y)$ 
are the Fourier transforms of $\phi_{0k}(y)$ and $\psi_{jk}(y)$
respectively. 
  
 Then % the process $X(t)$ can be represented as the series 
\begin{equation}
% W-РОЗКЛАД №1 !!!!!!
\label{p3eq-rozkl1}	
X(t)=\sum_{k \in \mathbb{Z}}\xi_{0k}a_{0k}(t)+
\sum_{j=0}^{\infty}\sum_{k \in \mathbb{Z}}\eta_{jk}b_{jk}(t),
\end{equation}
series $(\ref{p3eq-rozkl1})$ converges in $L_2(\Omega)$ for all
$\, t \in \mathbb{R}$, 
where $\xi_{0k},\, \eta_{jk}$  are  centered random variables   such that 
$$
\mathsf{E}\xi_{0k}\overline{\xi_{0l}}=\delta_{kl},\,\,\,
\mathsf{E}\eta_{jk}\overline{\eta_{lm}}=\delta_{jl}\delta_{km},\,\,\,
\mathsf{E}\xi_{0k}\overline{\eta_{nl}}=0.
$$ 
\end{thm}

\begin{nas} \cite{Koz-Turch-Roz}
\label{rozkl1-stats}
% === РОЗКЛАД №1 ДЛЯ СТАЦ. ПР-СІВ===
  Suppose that a centered second-order stationary process 
$X(t)$  has  the spectral density~$f(y)$, 
$\{ \phi_{0k}(x),\,\, \psi_{jk}(x),\,\,k \in \mathbb{Z},\,\,$
  $j=0,1, \ldots \}$ is a wavelet basis, $g(y)=\sqrt{f(y)}$. Then 
$X(t)$ can be represented as a  mean square convergent series 
$(\ref {p3eq-rozkl1})$ and 
\begin{equation}
\label{p3eq-rozkl1-a0k-stats}
a_{0k}(t)=\frac{1}{\sqrt{2\pi}} \int_{\mathbb{R}}g(y)
\exp\{-i y(t-k)\} \overline{\hat\phi(y)}dy,
\end{equation}  
\begin{equation}
\label{p3eq-rozkl1-bjk-stats}
b_{jk}(t)=\frac{1}{\sqrt{2\pi} \, 2^{j/2}} \int_{\mathbb{R}}
g(y) \exp\left\{-i y\left(t-\frac{k}{2^j}\right)\right\}
\overline{\hat\psi(y/2^j)}dy,
\end{equation}  
where the random variables $\xi_{0k},\, \eta_{jk}$ from  
$(\ref{p3eq-rozkl1})$ are such that
$$
\mathsf{E}\xi_{0k}=\mathsf{E}\eta_{jk}=0,
$$
$$
\mathsf{E}\xi_{0k}\overline{\xi_{0l}}=\delta_{kl},\,\,\,
\mathsf{E}\eta_{jk}\overline{\eta_{lm}}=\delta_{jl}\delta_{km},\,\,\,
\mathsf{E}\xi_{0k}\overline{\eta_{nl}}=0. 
$$ 
\end{nas}

\section{Simulation with given accuracy and reliability in $L_p([0,T])$ }

   By a stationary process we will always  mean a wide-sense stationary process below. 
   %Далі під стаціонарним процесом завжди будемо розуміти стаціонарний у
%широкому сенсі процес.   
%DEFI OF MODEL $\hat{X}$. DEFI OF SIMUL with given A \& R. 
%
%%
%We will  build a model which approximates a strictly sub-Gaussian random process with given accuracy and reliability in $C([0,T])$.
%% 

%%% {\bf   INTROD. PHRASES ?????!!!!! } 

%%%%%%% DEFIN OF Т. и Н.
\begin{ozn}
\label{def-model-hat-X}
  Suppose that a stationary random process $X=\{X(t),\,  t \in \mathbb{R}\}$
satisfies  the conditions of Corollary \ref{rozkl1-stats}.
 We call the following process a model of $X(t)$:   
\begin{equation}	
\label{pNR-model-hatX}
\hat {X}(t)=\sum\limits_{k=-(N_0-1)}^{N_0-1} \xi_{0k}a_{0k}(t)+
\sum\limits_{j=0}^{N-1}\sum\limits_{k=-(M_j-1)}^{M_j-1}\eta_{jk}b_{jk}(t),
\end{equation}
where $\xi_{0k}$, $\eta_{jk}$ are the random variables from the
expansion (\ref{p3eq-rozkl1}), $a_{0k}(t)$ and $b_{jk}(t)$ are calculated
using formulae (\ref{p3eq-rozkl1-a0k-stats}) and 
(\ref{p3eq-rozkl1-bjk-stats}), $N_0>1, N>1, M_j>1 \ (j=0, 1,  \ldots, N-1)$. 
\end{ozn}

  Numerical characteristics which describe the rate of approximation
of a process by its model are accuracy and reliability. 
  
\begin{ozn}
We say that a model $\hat{X}(t)$ approximates a process $X(t)$ with given  reliability 
$1-\delta$ ($0<\delta<1$) and  accuracy
$\varepsilon>0$ in $L_p([0,T]),\, p>0,$ if
\[
P \left\{\left( \int_0^T  |X(t) - \hat{X}(t)|^p \right)^{1/p} 
    > \varepsilon \right \} \leq \delta.
\]
\end{ozn}

\subsection{Simulation of $\mathbf{Y(t)=(X(t))^s}$}
 If $\hat{X}(t)$ is a model for a process 
$X(t)$ then a natural model for a process $Y(t)=F(X(t))$
is a ``plug-in'' model $\hat{Y}(t) = F(\hat{X}(t))$.  
So we will use $\hat{Y}(t) = (\hat{X}(t))^s$ as a model
 for $Y(t)=(X(t))^s$.  

\begin{thm}
\label{Main-Thm}
 Let  $Y(t)=(X(t))^s, t \in \mathbb{R}, s\in \mathbb{N},$ 
where  $X(t),  t \in \mathbb{R}, $ is a mean square 
continuous stationary strictly sub-Gaussian stochastic process 
which has spectral density $f(y)$, $g(y)=\sqrt{f(y)}$, 
$R(\tau)$ is the correlation function of $X(t)$, 
$\phi$ is a  $f$-wave\-let, $\psi$ is the corresponding  
$m$-wave\-let. Let the random variables $\xi_{0k}, \eta_{jk}$
in the expansion $(\ref{p3eq-rozkl1})$ of $X(t)$ be independent and 
strictly sub-Gaussian.  
Suppose that the following conditions hold:  there exist the derivatives 
$g'(y)$, $\hat \psi'(y)$, $\hat \phi'(y)$, 
$\vert \hat \psi(y)\vert <  C_1$, $\vert\hat \psi'(y)\vert < C_2$,
$f(y) \to 0$ as $|y| \to \infty$, $g(y)$ and $\hat \phi(y)$ are absolutely 
continuous,  
$$
\sup_{y \in \mathbb{R}} |\hat \phi(y)| < \infty, \quad 
\sup_{y \in \mathbb{R}} g(y) < \infty, 
$$
$$
\int_{\mathbb{R}} g(y)dy < \infty ,  \quad
 \int_{\mathbb{R}} \vert g'(y)\vert\vert y \vert dy < \infty ,
$$
$$
 \int_{\mathbb{R}} g(y)\vert y \vert dy < \infty \quad
 \int_{\mathbb{R}} \vert g'(y) \vert \vert \hat \phi(y) \vert dy <
 \infty ,
$$
$$
 \int_{\mathbb{R}} g(y) \vert \hat \phi'(y) \vert dy < \infty . 
$$ 
%%%%%%% 
%$\alpha_{0k}(t)$ and $\beta_{jk}(t)$ are given in (\ref{1.3})
%and (\ref{1.4}). If $k \not=0$ then for all $t \in \mathbb{R} $
%\begin{equation}
%\label{4.1}
%| \beta_{jk}(t)| \leq \frac{A+B|t|}{|k|2^{j/2}},
%\end{equation}
%where
Denote 
$$
A=\frac{C_2}{\sqrt{2\pi}} \int_{\mathbb{R}}(|g'(y)||y| + g(y))dy,
$$
$$
B=\frac{C_2}{\sqrt{2\pi}} \int_{\mathbb{R}} g(y)|y|dy, 
$$
%\begin{equation}
%\label{4.2}
%| \alpha_{0k}(t)| \leq \frac{A_1 + B_1|t|}{|k|},
%\end{equation}
%where
$$
A_1=\frac{1}{\sqrt{2\pi}}
\int_{\mathbb{R}} (|g'(y)||\hat\phi(y)| + g(y)|\hat\phi'(y)|)dy,
$$
$$
B_1=\frac{1}{\sqrt{2\pi}}
\int_{\mathbb{R}} g(y) |\hat\phi(y)|dy. 
$$
%%
% $$ 
% D=\frac{\textstyle{C}}{\textstyle{ \sqrt{2\pi} }}
  %\int_{\mathbb{R}} g(y)|y| dy.  
% $$
%%     {\bf ($C$ or $C_2$ in formula for $D$ ? $C_2$!!)(Yes, D=B !!!)}

  Let $p \geq 2, T>0, \varepsilon>0, 0<\delta<1$ and let the model 
$\hat{X}(t)$ of $X(t)$ be defined by $(\ref{pNR-model-hatX})$. Set 
$\hat{Y}(t) = (\hat{X}(t))^s$. 

  If  
\begin{equation}
\label{eq-N0-lb}
  N_0> \frac{6}{\delta_1}(A_1+B_1 T)^2 + 1,
\end{equation}
\begin{equation}
\label{eq-N-lb}
N>\max \left\{ 1 + \log_2 \left(\frac{72(A+B T)^2}
{5 \delta_1} \right), \,\,\,
1 + \log_8 \left( \frac{18B^2}{7 \delta_1} \right) \right \},
\end{equation}
\begin{equation}
\label{eq-Mj-lb}
M_j> 1 + \frac{12}{\delta_1}(A+B T)^2  \left( 1-\frac{1}{2^{N}} \right),
     \quad j= 0, 1,  \ldots,  N-1, 
\end{equation}
where 
$$
\delta_1 = \frac{\delta^{2/p} \varepsilon^2}{ (D^{*})^{2/p} }, 
$$
$$
D^{*}= 2^{(ps+3)/2} T p \sqrt{\Gamma(p)} s^{p-1/2} (R(0))^{p(s-1)/2}  
$$  
$$
\cdot \left[ (s-1)\Gamma(p(s-1)) +  
\sum_{k=1}^{s-2}\sqrt{k(s-1-k)\Gamma(2pk)\Gamma(2p(s-1-k))}\right]^{1/2},  
$$
then the model $\hat{Y}(t)$ approximates the process $Y(t)$ with given  
accuracy $\varepsilon$ and reliability $1-\delta$ in $L_p([0,T])$.  

%%%% SPRINGER !!!!
%%%% ==== НАСЛ-К О Т&Н ДЛЯ СТАЦ. СУБГАУССОВИХ, РОЗКЛ. №1 ======
%\begin{corollary}
%\label{nas-TN-subg-stats-roz1}
  %Suppose that a stationary sub-Gaussian random process 
 %$X= \{\mathbf{X}(t), \linebreak t \in \mathbb{R} \}$ satisfies the conditions of
%Corollary~\ref{pNR-cor-rozkl1-stats}, a $f$-wavelet $\phi$ and 
%the corresponding $m$-wavelet $\psi$ together with the process 
%$\mathbf{X}(t)$ satisfy the conditions of Lemma~\ref{p6lem-ocen-koef-rozkl1-stats}$;$ random variables
%$\{\xi_{0k},\,\,\eta_{jk},\,\,
%k\in  \mathbb{Z},\,\,j=0,1,\ldots\}$
%in expansion $(\ref{p3eq-rozkl1})$ of $\mathbf{X}(t)$ are independent
%and form a strictly sub-Gaussian family of random variables, 
%$p \geq 1;$ $A, A_1, B, B_1, C, \,g(y)$ are defined in 
%Lemma~\ref{p6lem-ocen-koef-rozkl1-stats}, $T>0$, 
%$$
%\delta_1=\min\left\{\frac{  \varepsilon^2 }
%{  2T^{2/p} \ln(2/\delta)},\,\,\frac{  \varepsilon^2 }
%{  p T^{2/p}} \right\}, \quad 
%D=\frac{\textstyle{C}}{\textstyle{ \sqrt{2\pi} }} \int_{\mathbb{R}} g(y)|y|\D y\; .  
%$$
   %If 
%$$
%N_0> \frac{6}{\delta_1}(A_1+B_1 T)^2 + 1,
%$$
%%
%$$
%N>\max \left\{ 1 + \log_2 \left(\frac{72(A+B T)^2}
%{5 \delta_1} \right), \,\,\,
%1 + \log_8 \left( \frac{18D^2}{7 \delta_1} \right) \right \},
%$$
%%
%$$
%M_j> 1 + \frac{12}{\delta_1}(A+B T)^2  \left( 1-\frac{1}{2^{N}} \right),
%$$
%then the model $\widehat{ \mathbf{X}}(t)$ defined 
%by~$(\ref{pNR-model-hatX})$  approximates the process
 %$\mathbf{X}(t)$ with reliability $1 - \delta$ and accuracy $\varepsilon$ in $L_p([0,T])$.   
%\end{corollary}
%%

\begin{proof}
  Denote 
	$\Delta X(t) = X(t) - \hat{X}(t)$, $\Delta Y(t) = Y(t) - \hat{Y}(t)$,
$$
\|\Delta Y\|_p = \left(\int_0^T |\Delta Y(t)|^p dt\right)^{1/p}. 
$$ 
  According to Lemma 4.1 from \cite{koz-turchyn-IJSMS}
the following inequalities hold under the conditions of the theorem: 
%%  ({\bf notations accord to IJSMS!!})
\begin{equation}
\label{eq-a0k-oc}
|a_{0k}(t)| \leq  \frac{A_1 + B_1|t|}{|k|},   \ k \neq 0, 
\end{equation}
\begin{equation}
\label{eq-bj0-oc}
|b_{j0}(t)| \leq  \frac{B}{2^{3j/2}}, 
%\frac{C_2}{\sqrt{2\pi}\, 2^{3j/2} } 
%\int_{\mathbb{R}} g(y)|y|dy = \frac{B}{2^{3j/2}},
\end{equation}
\begin{equation}
\label{eq-bjk-oc}
|b_{jk}(t)| \leq  \frac{A+B|t|}{|k|2^{j/2}},
 \ k \neq 0. 
\end{equation}
It follows from (\ref{eq-N0-lb})--(\ref{eq-Mj-lb}) and 
(\ref{eq-a0k-oc})--(\ref{eq-bjk-oc}) that 
$$
\sup_{t \in [0,T]} \mathsf{E}|\Delta X(t)|^2 \leq \delta_1. 
$$

  We will need the following inequality (see \cite{Rivasplata}): if $\xi$ 
is a sub-Gaussian random variable, then 
\begin{equation}
\label{oc-mom-subg}
  \mathsf{E} |\xi|^p \leq p\, 2^{p/2} (\tau(\xi))^p \, \Gamma(p/2), \ p>0.  
\end{equation}

   Let us estimate 
$\mathsf{E} \|\Delta {Y}  \|_p^2$.    
  Using % (\ref{oc-mom-subg}) and 
the Lyapunov inequality, we have:  
$$
\mathsf{E} \|\Delta {Y}  \|_p^2 =  
\mathsf{E}\left(\int_{0}^T  |Y(t) - \hat{Y}(t)|^p dt\right)^{2/p}  \leq
\left(\int_{0}^T \mathsf{E}|Y(t) - \hat{Y}(t)|^p dt\right)^{2/p}. 
$$	
% Let us estimate $\mathsf{E}|Y(t) - \hat{Y}(t)|^p$. 
Applying the Cauchy–Schwarz inequality we obtain: 
$$
\mathsf{E}|Y(t) - \hat{Y}(t)|^p = \mathsf{E}|(X(t))^s - (\hat{X}(t))^s|^p 
$$ 
$$
=  \mathsf{E}|X(t) - \hat{X}(t)|^p 
 \left| \sum_{k=0}^{s-1}  (X(t))^k (\hat{X}(t))^{s-1-k}\right|^p 
$$
$$
\leq  \left( \mathsf{E}|X(t) - \hat{X}(t)|^{2p}\right)^{1/2}
\left(\mathsf{E} \left| \sum_{k=0}^{s-1} 
 (X(t))^k (\hat{X}(t))^{s-1-k}\right|^{2p}\right)^{1/2} . 
$$
It follows from (\ref{oc-mom-subg}) that 
\begin{equation}
\label{mom-2p}
\left(\mathsf{E}|X(t) - \hat{X}(t)|^{2p}\right)^{1/2} \leq 
    (p\, 2^{p+1}\,\Gamma(p))^{1/2} (\sigma(\Delta X(t)))^p. 
\end{equation} 

Since 
\begin{equation}
\label{ineq-pow-mean}
 \mathsf{E}|\zeta_1 + \zeta_2 + \ldots + \zeta_s|^{2p} \leq 
  s^{2p-1} (\mathsf{E}|\zeta_1|^{2p} + 
	        \mathsf{E}|\zeta_2|^{2p} + \ldots + \mathsf{E}|\zeta_s|^{2p})					
\end{equation}
(an application of the power mean inequality) we have
\begin{equation}
\label{ineq-lsum1}
\mathsf{E} \left| \sum_{k=0}^{s-1} 
 (X(t))^k (\hat{X}(t))^{s-1-k}\right|^{2p} \leq  
s^{2p-1} \sum_{k=0}^{s-1} \mathsf{E}
 |X(t)|^{2pk} |\hat{X}(t)|^{2p(s-1-k)}. 
\end{equation}
But using (\ref{oc-mom-subg}) and Cauchy-Schwarz inequality we obtain    
\begin{equation}
\label{ineq-2mom}
\mathsf{E}  |X(t)|^{2pk} |\hat{X}(t)|^{2p(s-1-k)}
 \leq C_k^{*} (R(0))^{p(s-1)},  \quad  1\leq k\leq s-2,  
\end{equation}
where 
$$
  C_k^{*} = 
  4p\, 2^{p(s-1)} \sqrt{k(s-1-k)\Gamma(2pk)\Gamma(2p(s-1-k))}      
$$
(we used the inequality $\sigma(\hat{X}(t)) < \sigma(X(t))$). 

%Estimating the moments $\mathsf{E} |X(t)|^{2p(s-1)}$ and $\mathsf{E} |\hat{X}(t)|^{2p(s-1)}$
%using (\ref{oc-mom-subg}), we obtain 
It follows from (\ref{oc-mom-subg}) that   
\begin{equation}
\label{ineq-mom-X}
\mathsf{E}  |X(t)|^{2p(s-1)}  \leq  
2p(s-1) 2^{p(s-1)} \Gamma(p(s-1)) (R(0))^{p(s-1)}, 
\end{equation}
\begin{equation}
\label{ineq-mom-hatX}
\mathsf{E}  |\hat{X}(t)|^{2p(s-1)}  \leq  
2p(s-1) 2^{p(s-1)} \Gamma(p(s-1)) (R(0))^{p(s-1)}. 
\end{equation}
  Using (\ref{ineq-lsum1})--(\ref{ineq-mom-hatX}) we have  
$$	
\int_{0}^T \mathsf{E}|Y(t) - \hat{Y}(t)|^p dt \leq 
\delta_1^{p/2} D^{*} = \delta \varepsilon^p. 
$$
Applying the Markov inequality we get
$$
P\{\|\Delta Y\|_p > \varepsilon\} = 
P\{\|\Delta Y\|_p^p > \varepsilon^p \} \leq \delta.  
$$
So the theorem is proved. 
\end{proof}
\end{thm}

 %{\bf Do we estimate integr-s using grubuyu ots. cherez sup 
%ili integr-m ner-va ???!!!} 
%

 \begin{example}  
    A stationary centered Gaussian process  $X=\{ X(t),\, t \in \mathbb{R} \}$ 
 with spectral density $f(y)=\frac{\textstyle 1}{\textstyle(1+y^{2n})^2}, n \geq 2,$ 
and an arbitrary Daubechies wavelet satisfy the conditions of 
Theorem \ref{Main-Thm}. 
   %Прикладом процесу і вейвлету, які задовільняють умови теореми \ref{teor4.1},
%є центрований гауссовий стаціонарний процес $X=\{ X(t),\, t \in T \}$ з корреляційною функцією 
%$R(\tau)=e^{-\tau^2}$ і вейвлет Хаара (легко бачити, що для них 
%виконуються умови леми \ref{lemocmodnepr}).
  \end{example}

 %%% {\bf Remark. } Note that $\hat{Y}(t) \to {Y}(t)$  in probab??!

%\newpage
%\begin{center}
%{\Large {\bf SIMUL. OF $Z=X_1 X_2$}}
%\end{center}

\subsection{Simulation of $\mathbf{Z(t)=X_1(t) X_2(t)}$}

Let us now consider  a stochastic process $Z(t)$
which can be represented  as 
$$
  Z(t)=X_1(t) X_2(t), 
$$
where $X_1(t)$ and  $X_2(t)$ are independent  stationary strictly sub-Gaussian
stochastic processes which have spectral densities $f_1(y)$ and 
$f_2(y)$ correspondingly.   
Let $\phi_1, \psi_1$ and $\phi_2, \psi_2$ be two pairs of a $f$-wavelet and 
the corresponding $m$-wavelet. 

According to Corollary~\ref{rozkl1-stats}, 
the processes $X_1(t)$ and $X_2(t)$
can be expanded as 
\begin{equation}
% W-РОЗКЛАД №1 !!!!!!
\label{w1}	
X_1(t)=\sum_{k \in \mathbb{Z}}\xi_{0k}^{(1)}a_{0k}^{(1)}(t)+
\sum_{j=0}^{\infty}\sum_{k \in \mathbb{Z}}\eta_{jk}^{(1)}b_{jk}^{(1)}(t),
\end{equation}
\begin{equation}
% W-РОЗКЛАД №1 !!!!!!
\label{w2}	
X_2(t)=\sum_{k \in \mathbb{Z}}\xi_{0k}^{(2)}a_{0k}^{(2)}(t)+
\sum_{j=0}^{\infty}\sum_{k \in \mathbb{Z}}\eta_{jk}^{(2)}b_{jk}^{(2)}(t),
\end{equation}
where 
\begin{equation}
\label{ws-a0k-stats}
a_{0k}^{(s)}(t)=\frac{1}{\sqrt{2\pi}} \int_{\mathbb{R}} g_s(y)
\exp\{-i y(t-k)\} \overline{\hat\phi_s(y)}dy,
\end{equation}  
\begin{equation}
\label{bjk-s}
b_{jk}^{(s)}(t)=\frac{1}{\sqrt{2\pi} \, 2^{j/2}} \int_{\mathbb{R}}
g_s(y) \exp\left\{-i y\left(t-\frac{k}{2^j}\right)\right\}
\overline{\hat\psi_s(y/2^j)}dy,
\end{equation}  
$g_s(y)=\sqrt{f_s(y)},$ random variables $\xi_{0k}^{(s)}, \eta_{jl}^{(s)}$ are uncorrelated, 
$$
  \mathsf{E}|\xi_{0k}^{(s)}|^2=1,   \mathsf{E}|\eta_{jl}^{(s)}|^2=1,  
$$
$s=1, 2.$

  We will consider a ``plug-in'' model 
\begin{equation}
\label{mod-Zhat}
  \hat{Z}(t) = \hat{X}_1(t)\hat{X}_2(t)
\end{equation}	
for the process $Z(t)$, where 
$\hat{X}_1(t)$ and $\hat{X}_2(t)$ are models of type (\ref{pNR-model-hatX})
for $X_1(t)$ and $X_2(t)$ correspondingly, i.e. 
\begin{equation}
\label{mod-Xs} 
  \hat{X}_s(t)  = \sum\limits_{k=-(N_0^{(s)}-1)}^{N_0^{(s)}-1} \xi_{0k}^{(s)} a_{0k}^{(s)}(t)+
\sum\limits_{j=0}^{N^{(s)}-1}\sum\limits_{k=-(M_j^{(s)}-1)}^{M_j^{(s)}-1}\eta_{jk}^{(s)} b_{jk}^{(s)}(t),
\end{equation}	
$N_0^{(s)}>1, N^{(s)}>1, M_j^{(s)}>1 \ (j=0, 1,  \ldots, N^{(s)}-1), \ s=1, 2.$
%
%\begin{equation}
%\label{mod-X2}
  %\hat{X}_2(t)  =\ldots  
%\end{equation}	

% \baselineskip=20pt 

\begin{thm}
\label{AR-X1X2}
 Let  $Z(t)=X_1(t)X_2(t), t \in \mathbb{R}$,  
where  $X_1(t), X_2(t), t \in \mathbb{R},$ are mean square 
continuous stationary strictly sub-Gaussian stochastic processes  
which have spectral densities $f_s(y)$, $g_s(y)=\sqrt{f_s(y)}$, $R_s(\tau)$ is 
the correlation function of $X_s(t) (s=1,2)$,   
$\phi_1, \psi_1$ and  $\phi_2, \psi_2$ are two pairs of a 
$f$-wave\-let and the corresponding  $m$-wave\-let. 
Let the random variables $\xi_{0k}^{(s)}, \eta_{jk}^{(s)} (s=1,2)$
in expansions $(\ref{w1}), (\ref{w2})$ of $X_1(t)$ and $X_2(t)$
be independent and strictly sub-Gaussian.  
% {\bf  CHECK CAREFULLY !!} 
Suppose that the following conditions hold:  there exist derivatives 
$g_s'(y)$, $\hat{\psi}_s'(y)$, $\hat{\phi}_s'(y)$, 
$\vert \hat{\psi}_s(y)\vert <  C_1^{(s)}$,
 $\vert\hat{\psi}_s'(y)\vert < C_2^{(s)}$,
$f_s(y) \to 0$ as $|y| \to \infty$, $g_s(y)$ and $\hat{\phi}_s(y)$ are absolutely 
continuous,   
$$
\sup_{y \in \mathbb{R}} |\hat{\phi}_s(y)| < \infty, \quad 
\sup_{y \in \mathbb{R}} g_s(y) < \infty, 
$$
$$
\int_{\mathbb{R}} g_s(y)dy < \infty ,  \quad
 \int_{\mathbb{R}} \vert g_s'(y)\vert\vert y \vert dy < \infty ,
$$
$$
 \int_{\mathbb{R}} g_s(y)\vert y \vert dy < \infty \quad
 \int_{\mathbb{R}} \vert g_s'(y) \vert \vert \hat{\phi}_s(y) \vert dy <
 \infty ,
$$
$$
 \int_{\mathbb{R}} g_s(y) \vert \hat{\phi}_s'(y) \vert dy < \infty,  
$$ 
$s=1,2.$
%%%%%%% 
%$\alpha_{0k}(t)$ and $\beta_{jk}(t)$ are given in (\ref{1.3})
%and (\ref{1.4}). If $k \not=0$ then for all $t \in \mathbb{R} $
%\begin{equation}
%\label{4.1}
%| \beta_{jk}(t)| \leq \frac{A+B|t|}{|k|2^{j/2}},
%\end{equation}
%where
Denote 
$$
A^{(s)}=\frac{C_2^{(s)}}{\sqrt{2\pi}} \int_{\mathbb{R}}(|g_s'(y)||y| + g_s(y))dy,
$$
$$
B^{(s)}=\frac{C_2^{(s)}}{\sqrt{2\pi}} \int_{\mathbb{R}} g_s(y)|y|dy, 
$$
%\begin{equation}
%\label{4.2}
%| \alpha_{0k}(t)| \leq \frac{A_1 + B_1|t|}{|k|},
%\end{equation}
%where
$$
A_1^{(s)}=\frac{1}{\sqrt{2\pi}}
\int_{\mathbb{R}} (|g_s'(y)||\hat{\phi}_s(y)| + g_s(y)|\hat{\phi}_s'(y)|)dy,
$$
$$
B_1^{(s)}=\frac{1}{\sqrt{2\pi}}
\int_{\mathbb{R}} g_s(y) |\hat{\phi}_s(y)|dy,  
$$
$s=1,2.$
%%
% $$ 
% D=\frac{\textstyle{C}}{\textstyle{ \sqrt{2\pi} }}
  %\int_{\mathbb{R}} g(y)|y| dy.  
% $$
%%     {\bf ($C$ or $C_2$ in formula for $D$ ? $C_2$!!)(Yes, D=B !!!)}

  Let $p \geq 1, T>0, \varepsilon>0, 0<\delta<1$ 
% ({\bf check cond-s on $p$!})
and let the models 
$\hat{Z}(t), \hat{X}_s(t)$ of 
$Z(t), X_s(t)$ $(s=1, 2)$ be defined by 
$(\ref{mod-Zhat}), (\ref{mod-Xs})$ correspondingly. 
  
	If 
\begin{equation}
\label{eq-N0s-lb}
  N_0^{(s)} > \frac{6}{\delta_s^{*}}(A_1^{(s)}+B_1^{(s)} T)^2 + 1,
\end{equation}
\begin{equation}
\label{eq-Ns-lb}
N^{(s)}>\max \left\{ 1 + \log_2 \left(\frac{72(A^{(s)}+B^{(s)} T)^2}
{5 \delta_s^{*}} \right), \,\,\,
1 + \log_8 \left( \frac{18(B^{(s)})^2}{7 \delta_s^{*}} \right) \right \},
\end{equation}
\begin{equation}
\label{eq-Mjs-lb}
M_j^{(s)}> 1 + \frac{12}{\delta_s^{*}}(A^{(s)}+B^{(s)} T)^2  \left( 1-\frac{1}{2^{N^{(s)}}} \right),
     \quad j= 0, 1,  \ldots,  N^{(s)}-1, 
\end{equation}
$s=1,2,$  where 
$$
\delta_1^{*} =  \hat{\delta}/R_2(0), \  \delta_2^{*} =  \hat{\delta}/R_1(0), 
$$
$$
 \hat{\delta} = \frac{ \delta^{2/p}\varepsilon^2 }{ (2^{2p+1}p\Gamma(p)T)^{2/p} }, 
$$
 then the model $\hat{Z}(t)$  defined by $(\ref{mod-Zhat})$ approximates the process $Z(t)$ with given  
accuracy $\varepsilon$ and reliability $1-\delta$ in $L_p([0,T])$.  
\begin{proof}

 Denote $\Delta X_s(t)=X_s(t) - \hat{X}_s(t),$
$$
\|\Delta X_s\|_p = \left(\int_0^T |\Delta X_s(t)|^p dt\right)^{1/p},  
$$ 
$s=1,2.$    
 Let us estimate 
	$$
	  \mathsf{E} \|X_1X_2 -\hat{X_1}\hat{X_2}\|_p. 
	$$
	
Applying  (\ref{oc-mom-subg}), the Cauchy-Schwarz inequality and the power mean inequality we get 
$$
   \mathsf{E}|X_1(t)X_2(t) -\hat{X_1}(t)\hat{X_2}(t)|^p 
$$
%$$	
	%\leq 
   %\mathsf{E}(|(X_1(t)X_2(t) - \hat{X_1}(t)X_2(t)|  + 
		          %|\hat{X_1}(t)X_2(t) - \hat{X_1}(t)\hat{X_2}(t)|)^2 
%$$
$$
\leq 	 2^{p-1} (\mathsf{E}|X_2(t)(X_1(t) - \hat{X_1}(t))|^p   + 
           \mathsf{E}|\hat{X_1}(t)(X_2(t) - \hat{X_2}(t))|^p)
$$					
$$
\leq 	 2^{p-1} \Big((\mathsf{E}|X_2(t)|^{2p})^{1/2} (\mathsf{E}|X_1(t) - \hat{X_1}(t)|^{2p})^{1/2}  
$$
$$
 +  (\mathsf{E}|\hat{X_1}(t)|^{2p})^{1/2} (\mathsf{E} |X_2(t) - \hat{X_2}(t)|^{2p})^{1/2}\Big)
$$
$$
\leq 2^{p-1} \Big( C_{R}(p) \sigma^p(X_2(t)) \sigma^p(X_1(t) - \hat{X_1}(t))
%$$
%\begin{equation*}
%\label{eq-MSE-X1X2}
 +  C_{R}(p) \sigma^p(X_1(t)) \sigma^p(X_2(t) - \hat{X_2}(t)) \Big),    
$$
%\end{equation*}
where 
$$
 C_{R}(p) = 2p 2^p \Gamma(p). 
$$

It follows from 
(\ref{eq-N0s-lb})--(\ref{eq-Mjs-lb}) that 
$$
 \sup_{t \in [0;T]} \sigma^2 (\Delta X_s(t)) \leq \delta_s^*
$$
$s=1, 2.$  Now we see that  
$$	
\int_{0}^T \mathsf{E}|Z(t) - \hat{Z}(t)|^p dt 
$$
$$
\leq 2^{p-1} C_{R}(p) T 
\Big(
 (R_2(0))^{p/2} \big( \sup_{t \in [0,T]} 
             \sigma^2(\Delta X_1(t))\big)^{p/2}
+ (R_1(0))^{p/2}  \big( \sup_{t \in [0,T]}  \sigma^2(\Delta X_2(t))\big)^{p/2}
 \Big)
$$
$$
\leq      \delta \varepsilon^p 
$$ 
and,  using the Markov inequality, we obtain % SYN: obtain/get 
$$
P\{\|Z - \hat{Z}\|_p > \varepsilon\}  \leq \delta.  
$$
\end{proof}
\end{thm}

\begin{example}  
% {\bf Example.}
Let us consider a process $Z(t)=X_1(t) X_2(t)$, 
where $X_1(t)$ and $X_2(t)$ are independent centered stationary Gaussian stochastic
processes which have spectral densities 
$$
f_1(y) = \frac{1}{(1+y^2)^{2n} },  \  n \geq 2,
$$ 
and 
%% modulated 
$$
f_2(y) = \left(\frac{1}{(1+(y-a)^2)^m} + \frac{1}{(1+(y+a)^2)^m}\right)^2, \  m \geq 2, 
$$ 
correspondingly.  Let us take as $\phi_1, \phi_2$ 
and $\psi_1, \psi_2$  correspondingly two Daubechies $f$-wavelets and $m$-wavelets of any order.
 These two pairs of a process and the corresponding wavelet satisfy the conditions 
of the theorem.       
%spectr. dens. $f_1(y) = 1/(1+y^2)^n, n \geq ...,$ and (modulated dens-ty)
%$$
%f_2(y) = \frac{1}{(1+(y-a)^2)^m} + \frac{1}{(1+(y+a)^2)^m}, \  m \geq ...
%$$ 
\end{example}

\section{Conclusions}
   
	 We built a wavelet-based model for simulation of a process which is an 
	integer power of a sub-Gaussian process. A wavelet-based model was also built 
for a process $Z(t)$ which can be represented as $Z(t)=X_1(t) X_2(t)$, where 
$X_1(t)$ and $X_2(t)$ are stationary strictly sub-Gaussian processes.  
	
	We proved theorems about simulation of stochastic processes by the above-mentioned models 
  with given accuracy and reliability in $L_p([0;T])$.

   The author expresses gratitude to professor Yury Kozachenko for valuable 
  discussions.

%% =====NEW=======

%Ayache, A., Linde, W. (2008). 
%\textit{Approximation of Gaussian random fields: 
%general results and optimal wavelet representation
%of the L?vy fractional motion}
%Journal of Theoretical Probability,
 %21(1), 69--96.

%\bibitem{Ayache-LFSS}
%A. Ayache, F. Roueff, Y. Xiao,
%\textit{Linear fractional stable sheets: Wavelet expansion and
%sample path properties},
%Stochastic Processes and their Applications
 %119 (2009), 
%1168--1197. 

%\bibitem{paper}
%A. U. Thor,
%\textit{Article title},
%Journal Name
%\textbf{17}
%(2000),
%no. 1,
%111--113.

\end{document}